\documentclass{conm-p-l}
\usepackage{amssymb}
\usepackage{epic,eepic}

\newcommand\eerf[1] {Eq.\ \nolinebreak (#1)}
\def\ch            {{\chi}}
\def\a             {\alpha}

\def\be            {\begin{equation}}

\def\bbZ           {\mathbb{Z}}

\def\cA            {{\mathcal{A}}}
\def\cB            {{\mathcal{B}}}
\def\cC            {{\mathcal{C}}}

\def\cH            {{\mathcal{H}}}

\def\cP            {{\mathcal{P}}}

\def\cS            {{\mathcal{S}}}
\def\cT            {{\mathcal{T}}}

\def\ee            {\end{equation}}

\newcommand\erf[1] {Eq.\ \nolinebreak (\ref{#1})}

\def\ext           {{\mathrm{ext}}}
\newcommand\fig[1] {Fig.\ \nolinebreak \ref{#1}}

\def\I             {{\mathrm{i}}}

\def\la            {\lambda}
\def\lan           {\langle}

\def\ran           {\rangle}
\def\rmA           {{\mathrm{A}}}
\def\rmD           {{\mathrm{D}}}
\def\rmE           {{\mathrm{E}}}

\def\SLZ           {{\mathit{SL}}(2;\bbZ)}

\def\SUn           {{\mathit{SU}}(n)}

\def\SUz           {{\mathit{SU}}(2)}

\def\tr            {{\mathrm{tr}}}
\def\Tr            {{\mathrm{Tr}}}

\newtheorem{theorem}{Theorem}[section]

\theoremstyle{definition}

\theoremstyle{remark}

\numberwithin{equation}{section}

\begin{document}

\title{Modular Invariants and their Fusion Rules}

\author{David E Evans}
\address {School of Mathematics, University of Wales Cardiff,
Senghennydd Road, Cardiff CF24 4YH, Wales, United Kingdom}
\email{EvansDE@cf.ac.uk}
\author{Paulo R Pinto}
\address{Department of Mathematics, Inst.\ Super.\ T\'ecnico, Av.\
Rovisco Pais, 1049-001 Lisboa, Portugal}
\email{ppinto@math.ist.utl.pt}
\thanks{The second author was partially supported by FCT (Portugal) grant BD/9704/96 and CAMGSD-IST}

\date{1 November, 2002 and, in revised form, *******.}

\dedicatory{This paper is dedicated to Robert T. Powers on the occasion of
his sixtieth birthday}


\begin{abstract}
The subfactor approach to modular invariants gives insight
into the fusion rule structure of the modular invariants.
\end{abstract}

\maketitle

\section{Introduction and Background}

We are going to use the inter-relations between quantum and
classical dynamical systems and to use tools from non-commutative
geometry or non-commutative operator algebras to understand
phenomena in classical statistical mechanics.

A good starting point is the two dimensional Ising model and its
C*-treatment \cite{AE,EL,EK}. The classical model is set in the
configuration space $\cP$ = $\{\pm\}^{{\mathbb{Z}^2}}$ of
distributions of $+$ and $-$ on the two dimensional lattice
${\mathbb{Z}^2}$, and the nearest neighbour Hamiltonian

\[
H=-\sum_{\alpha, \beta
~nn}J\sigma^\alpha\sigma^\beta
\]

\noindent where the sum is over nearest neighbours $\alpha$ and
$\beta$ in $\mathbb{Z}^2$, and $\sigma  = (\sigma^{\alpha})$ in
$\cP$. The transfer matrix method allows us to study the
classical model set in $C(P)
=\bigotimes_{\mathbb{Z}^2}{\mathbb{C}}^2$
 and its equilibrium states, characterized by say the
Dobrushin-Lanford-Ruelle equations or a variational
principle by a quantum system of noncommuting observables
$\mathcal{A}=\bigotimes_{\mathbb{Z}}M_2$
 in one dimension with dynamics
$\alpha_t = T^{it}( - )T^{-it}$ and associated equilibrium states
or more precisely its ground states. Associated to an equilibrium
state $\mu$ at inverse temperature $\beta$ for the Ising model is
a ground state $\varphi_{\mu}$ on the Pauli algebra $\cA$ and to
each local observable $F$ in $C(\cP)$, a quantum observable
$F_{\beta}$ depending only on the temperature such that we can
describe the classical correlation values in terms of quantum
ones: $\mu(F) = \varphi_{\mu}(F_\beta)$.
The key element in this reduction is to identify the classical
partition $Z$ function in the denominator of the classical
expectation value
\begin{equation}\label{partition}
\mu(F)=\sum_{\sigma} F(\sigma) \exp(-\beta
H(\sigma))/\sum_{\sigma} \exp(-\beta H(\sigma))
\end{equation}
\noindent in quantum terms. We write for horizontal periodic
boundary conditions:

${}$

 $Z\  =\  \sum_{\sigma} \exp(-\beta H(\sigma))$

\begin{table}[hhhhhhhh]
$$
\hskip -2.0cm =\quad
\begin{array}{|c|c|c|c|c|c|}
\hline {}&{}&{}&{}&{}&{}\\ \hline &&&&&\\ \hline
{}&{}&{}&{}&{}&{}\\ \hline &&&&&\\ \hline {}&{}&{}&{}&{}&{}\\
\hline
\end{array}
\quad =\sum_\sigma\quad
\begin{array}{|c|}
\hline  \\ \hline   \\ \hline \\ \hline
\\ \hline \\ \hline
\end{array}
\quad
\begin{array}{|c|}
\hline  \\ \hline   \\ \hline \\ \hline
\\ \hline \\ \hline
\end{array}
\quad
\begin{array}{|c|}
\hline  \\ \hline   \\ \hline \\ \hline
\\ \hline \\ \hline
\end{array}
\quad
\begin{array}{|c|}
\hline  \\ \hline   \\ \hline \\ \hline
\\ \hline \\ \hline
\end{array}
\quad
\begin{array}{|c|}
\hline  \\ \hline   \\ \hline \\ \hline
\\ \hline \\ \hline
\end{array}
\quad
\begin{array}{|c|}
\hline  \\ \hline   \\ \hline \\ \hline
\\ \hline \\ \hline
\end{array}
$$
\end{table}
\[ =  \sum_\sigma T{_{\xi_1\xi_2}}T{_{\xi_2\xi_3}}\ldots
T_{\xi_N\xi_1}={\rm trace }\ T \dots T = {\rm trace}\ T^N =
{\rm trace}\ {\rm e}^{-N\cH}. \]


The partition function of a square of size $M$ by $N$ decomposes
as the exponential of a sum factorising into a product of the
partition functions of columns summed over column configurations
and hence can be understood as a matrix product, computed in the
trace as we have imposed periodic boundary conditions. Taking the
scaling limit at criticality, we obtain a field theory, and for a
torus model with periodic boundary conditions in addition
vertically, we have
\begin{equation}
Z(\tau) =
\tr(q^{(L_0-c/24)}\bar{q}^{(\bar{L}_0-c/24)}) \,.
\label{parf1}
\end{equation}

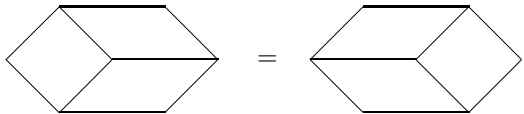
\begin{figure}[tb]
\begin{center}
\begin{picture}(215,60)
\path(10,30)(30,50)(50,30)(30,10)(10,30)
\put(30,50){\line(1,0){40}}
\put(30,10){\line(1,0){40}}
\put(50,30){\line(1,0){40}}
\put(70,50){\line(1,-1){20}}
\put(70,10){\line(1,1){20}}
\put(105,28){$=$}
\multiput(125,30)(40,0){2}{\line(1,1){20}}
\multiput(125,30)(40,0){2}{\line(1,-1){20}}
\put(125,30){\line(1,0){40}}
\multiput(145,10)(0,40){2}{\line(1,0){40}}
\put(185,50){\line(1,-1){20}}
\put(185,10){\line(1,1){20}}
\end{picture}
\end{center}
\caption{Yang-Baxter equation}
\label{figYBE}
\end{figure}

\noindent Here $L_0$ and $\bar L_0$ are commuting Hamiltonians, usually part of
commuting
Virasoro algebras as for the Ising model here with associated
central charge $c$, or multiplier of the projective representation
$L_m$ of the vector fields $l_m$ = $-z^{m+1}d/dz$ on the circle.
Much of this structure can be understood or is present in the
statistical mechanical model itself, and one of our aims is to
lift structures or understanding of the conformal field theories
at criticality back to the original statistical mechanical models.
For example, the central charge $c$ itself can be seen from the
asymptotics of the partition function $Z$:
  when $1 << N << M$, then
$Z \approx \exp (-NMf + {M \pi c} / N 6)$, where $f$ is the free
energy ${\rm log}(Z)/MN$ after $N,M \rightarrow \infty$
\cite{cardy2}.
 To generalise the Ising model we need some integrability of
solutions of the Yang-Baxter equation as in  Fig.\ \ref{figYBE}.
for local Boltzmann weights.
Such solutions are naturally provided by
Hecke algebras or quantum groups associated to $\SUn_k$ particularly
at roots of unity - so that $\SUz_2$, the Ising model is the first in a
double series of examples.

The Boltzmann weights lie in
\[
(\bigotimes M_n)^{\SUn_q} = \pi({\rm Hecke}).
\]
The justification of the term $\SUn$ models is as follows. By
Weyl duality, the representation of the permutation group in
$\bigotimes M_n$ is the fixed point algebra of the product action
of $\SUn$. Deforming this, there is a representation of the Hecke
algebra in $\bigotimes M_n$ whose commutant is a representation
of a deformation of $\SUn$, the quantum group $\SUn_q$ \cite{Ji}.
 The braid
relations are then precisely the YBE at criticality. The
representations which appear here are labelled by Young tableaux
of at most $n-1$ rows but a further constraint implying
rationality or finiteness of the representation labels occur at
roots of unity $q = e^{\pi i/(n+k)}$ when only labels with at
most $k$ columns appear. The $\SUz_k$ models are constructed by
distributing edges of the Dynkin diagram $\rmA_{k+1}$ on square
lattice. The first non trivial example here is described by
$\rmA_{2}$, the two extreme vertices representing our previous
symbols $\{\pm \}$ and the internal symbol is a dummy variable.
The labels thus have two meanings coming from the representation
theory of $\SUn$ of symmetric group (and their deformations) and
which point of view we want to emphasise may depend on whether we
are more interested in the statistical mechanics or the conformal
field theory picture:

\begin{center}

$\begin{array}{|c|c|}\hline {}&{}\\ \hbox{statistical
mechanics}&\hbox{conformal field theory} \\{}&{}
\\\hline
{}&{}\\
\hbox{Hecke algebras}&SU(N)\hbox{ loop groups}
\\ {}&{}
\\ \hline
{}&{}\\ \la \hbox{ symmetric or braided}\atop {} \hbox{group
representation}& \la \hbox{ positive energy}\atop {}
\hbox{representation of $SU(N)$}
\\ {}&{}
\\ \hline
{}&{}\\
\hbox{II$_1$ bimodules}\atop \hbox{Jones-Wenzl
subfactors}&\hbox{III$_1$ sectors}\atop \hbox{Jones-Wassermann
subfactors}
\\
{}&{}
\\ \hline  {}&{}\\ \hbox{Ocneanu paragroup}\atop \hbox{and
connections}&\hbox{Longo Q-systems and}\atop\hbox{canonical
endomorphism}\\ {}&{}
\\ \hline
\end{array}$
\end{center}

\section{Modular Invariants}

The link between the two frameworks are deformed Weyl duality,
Popa's classification of type III inclusion in terms of those of
type II and of course the partition function $Z$. Decomposing
\erf{parf1} according to our underlying loop group symmetries we
have:
\[ Z(\tau) = \sum\nolimits_{\la,\mu} Z_{\la,\mu}
\chi_\la(\tau) \chi_\mu(\tau)^* \,,\] \noindent where
$\chi_\la(\tau) = \tr_{H_\la} \exp (2\pi\I \tau L_0)$ is the
character in the positive energy representation $\lambda$ and $q
= e^{2 \pi i\tau}$. Modular invariance under reparameterisation
of the torus then becomes $Z(\tau) = Z((a\tau + b)/(c \tau + d))$
for $\left({a\atop c}{b \atop d}\right)$ in $\SLZ$ which is
generated by $\cS=\left({0\atop 1}{-1\atop 0}\right)$ and
$\cT=\left({1\atop 0}{1\atop 1}\right)$. Since $\SLZ$ acts
linearly on characters, there is a representation of $\SLZ$ by
$\cS \mapsto S=[S_{\lambda,\mu}]$, and $\cT \mapsto
T=[T_{\lambda,\mu}]$,  where
\[
\chi_\la(-1/\tau)= \sum\nolimits_\mu S_{\la,\mu} \chi_\mu(\tau) \,,
\qquad \chi_\la(\tau+1) = \sum\nolimits_\mu T_{\la,\mu} \chi_\mu(\tau) \,.
\]
Then modular invariance is
\begin{equation}\label{help}
ZS = SZ, ZT = TZ.
\end{equation}

\noindent Since the coefficients $Z_{\lambda \mu}$
appeared as multiplicities,
\begin{equation}\label{integral}
 Z_{\lambda \mu}  \in 0,1,2, \dots,
\end{equation}
 and
usually through uniqueness of the vacuum we have
\begin{equation}\label{normalization}
 Z_{0 0} =1.
\end{equation}
A modular invariant is then a matrix $Z = [Z_{\lambda \mu}]$
satisfying the above constraints \erf{help}, \erf{integral}  and
\erf{normalization}, although we will see that much interesting
structure is uncovered if we relax the normalization condition
$Z_{0 0} =1$. Usually as in the $\SUn_k$ case, the modular data
is such that $S$ is symmetric and $S_{0 \lambda} >0$.

In this case \cite{BER1,G1}
\begin{equation}\label{Zbound}
0 \le Z_{\la,\mu} \le d_\la d_\mu
\end{equation}

\noindent where $d_\lambda = S_{0,\lambda}/S_{00}$, and thus there are at most
finitely many normalized modular invariants.

There always exists the trivial modular invariant $Z = 1$. For
$\SUz$ there is the celebrated A-D-E classification \cite{CIZ2},e.g.
$\SUz_{16}$ have three modular invariants labelled by the Dynkin
diagrams  $ \rmA_{17}, \rmD_{10}, \rmE_{7}$ with corresponding
invariants $Z_{\rmA_{17}} = 1$, $Z_{\rmD_{10}}$ is an orbifold
obtained by folding the $\rmA_{17}$-diagram which has a fixed
point, and $Z_{\rmE_7}$ is a twist of the orbifold invariant.
These are the normalized modular invariants, and if one relaxes
the normalization condition, $Z_{00} =1$, to consider for example
products of normalized modular invariants then an interesting
fusion rule algebra appears as in Table \ref{ff}.

\unitlength 1.4mm \thinlines
\begin{table}[htb]
\begin{center}
\begin{tabular}{|c|c c c |} \hline
&$Z_{\rmA_{17}}$&$Z_{\rmD_{10}}$&$Z_{\rmE_{7}}$\\    \hline\hline
$Z_{\rmA_{17}}$& $Z_{\rmA_{17}}$&$Z_{\rmD_{10}}$&$Z_{\rmE_{7}}$\cr
$Z_{\rmD_{10}}$&$Z_{\rmD_{10}}$&$2Z_{\rmD_{10}}$&$2Z_{\rmE_{7}}$\cr
$Z_{\rmE_{7}}$&$Z_{\rmE_{7}}$&$2Z_{\rmE_{7}}$
&$Z_{\rmD_{10}}+Z_{\rmE_{7}}$\cr \hline
\end{tabular}
\end{center}
\caption{Fusion of $\SUz_{16}$ modular invariants} 
\label{ff}
\end{table}

It is then natural to ask whether the (unnormalized) modular
invariants form a fusion rule algebra generated by the normalized
modular invariants and get a better understanding of these fusion
rule algebras. The subfactor approach can assist us in this \cite{E1,EP}.

\section{Moore Seiberg dilation}

Here we will have a type III factor $N$, endowed with a system
$\cA$ of endomorphisms which are taken to be braided and so yield
a modular data.
As we have seen we always have the trivial modular invariant
\begin{equation}
Z= \sum\nolimits_\la |\chi_\la|^2,
\end{equation}
where we now interprete $\chi_{\lambda}$ as a formal character
${\rm tr}q^{L_0 - c/24}$
even when the Hamiltonian $L_0$ may not exist. More generally we
may introduce twists or permutations of the fusion rules preserving
$S$ and $T$ and the vacuum $0$ so that we should also consider
\begin{equation}\label{twist}
Z = \sum\nolimits_\la \chi_\la \chi_{\vartheta(\la)}^* .
\end{equation}
In some sense every modular invariant can be dilated or brought
to this form. Suppose we can extend the system in the following
sense. We have a subfactor $N \subset M$ with a  system $\cB$ of
endomorphisms of $M$. One can emphasise the extension aspect of
moving from endomorphisms of $N$ to those of $M$
($\alpha$-induction \cite{LR,BE2}). One needs the system $\cA$ to
be braided and the sector of $\bar\iota\iota$, the dual canonical
endomorphism to lie in $\Sigma(\cA)$, the set of finite integral
sums of endomorphism of $\cA$.

Under these conditions we emphasise the restriction aspect of this
extension. The trivial or twisted modular invariants \erf{twist}
for the $\cB$-system can
restricted to $\cA$, written formally in terms of characters
as
$\chi_{\tau} = \sum_{\lambda}b_{\tau,\lambda} \,\chi_\lambda $,
with branching coefficients $b_{\tau,\lambda}$. Then  restricting  the
diagonal modular invariant
\begin{equation}
Z^{\ext}=\sum\nolimits_{\tau\in\cB} |\chi_{\tau}|^2,
\end{equation}
\noindent to the original system  we have:

\begin{equation}
\sum_{\tau\in
\mathcal{B}}|\chi_\tau|^2 = \sum_{\tau\in
\mathcal{B}}|\sum_{\lambda\in
\mathcal{A}}b_{\tau\lambda}\chi_\lambda|^2,
\end{equation}


\noindent with mass matrix

\begin{equation}
\label{ext}
Z_{\la,\mu}=\sum\nolimits_\tau
b_{\tau,\la} b_{\tau,\mu}.
\end{equation}

\noindent These invariants are called type I (or more precisely
the inclusion $N \subset M$ describing this modular invariant is
type I) and are necessarily symmetric $Z_{\lambda\mu} =
Z_{\mu\lambda}$. In the presence of non-trivial twist $\vartheta$
of the $\cB$ system, we have type II invariants

\begin{equation}
\label{extaut}
Z_{\la,\mu}=\sum\nolimits_\tau
b_{\tau,\la} b_{\vartheta(\tau),\mu}.
\end{equation}

\noindent which have symmetric vacuum coupling $Z_{0\lambda} =
Z_{\lambda 0}$.

A type III phenomena accurs, due to some underlying heterotic 
structure which results in needing different labellings $\cB_+$ and $\cB_-$ 
on left and right extended systems $\cA\subset \cB_\pm$. In the subfactor 
framework, this is found in two intermediat
e subfactors $N\subset M_\pm\subset M$, where $M_\pm$ carry systems of 
endomorphisms $\cB_\pm$. Not only can this situation be found with modular 
invariants with non-symmetric vacuum coupling as in the orthogonal 
at loop groups low levels \cite{BE4}, but also with quantum doubles of
finite groups \cite{EP}.
Turning to the extension point of view $N\subset M$, when 
$\cA\subset \hbox{End}(N)$ is non-degenerately braided and the dual 
canonical endomorphism $\theta$ lies in $\Sigma(\cA)$, we form $\cC^\pm$ to 
be the systems of endomorphisms on $M$ from 
the irreducible components of the sectors of $\{\a_\la^\pm: \la\in\cA\}$. 
These generate the full system $\cC=\cC^+\vee\cC^-$, not necessarily even 
commutative, but $\cB=\cC^0=\cC^+\cap\cC^-$ is not only commutative but 
non-degenerately braided \cite{BE4}.


\section{Izumi quantum $E_6$ model}

Non-degenerately braided systems and the corresponding modular data can 
be obtained in the operator algebra setting in at least two ways. One is
to take the Jones-Wassermann loop group examples in algebraic quantum
field theory analysed by Wasserman \cite{W} and students e.g.
Laredo-Toledano
\cite{Tol}. The other is to take the quantum double of a system of 
endomorphisms - which may not even be commutative let alone braided
or may be braided but the braiding has some degeneracy. Here we focus
on one such example the quantum double of the even $E_6$ system. (As 
explained in \cite[Section 5]{BEK3}, this is more natural than taking 
the double of the entire $E_6$ system as the Longo-Rehren 
construction then gives the double of the $E_6$ subfactor). 
The $E_6$ system can be realised
through a conformal embedding of 
$\mathit{SU}(2)_{10}\subset \mathit{SO}(5)_1$ \cite{X2, BE2, BE3}. 
The  modular data obtained from the Longo-Rehren
inclusion of the even system of  $E_6$ is 
explicitly written down as follows in \cite[Page 648]{iz3}
where $d=1+\sqrt{3}, \la=2+d^2, i=\sqrt{-1}.$ The $S$ and $T$
matrices are

{\scriptsize
\begin{eqnarray}
S&=&\frac{1}{\la} \begin{pmatrix} 1&1&1+d&1+d&2+d&d&d&d&d&d\cr
1&1&1+d&1+d&-2-d&d&d&-d&-d&-d\cr 1+d&1+d&1&1&2+d&-d&-d&-d&-d&-d\cr
1+d&1+d&1&1&-2-d&-d&-d&d&d&d\cr 2+d&-2-d&2+d&-2-d&0&0&0&0&0&0\cr
d&d&-d&-d&0&d&d&-d&-d&2d\cr d&d&-d&-d&0&d&d&-d&-d&-2d\cr
d&-d&-d&d&0&-d&d&-2i-di&2i+di&0\cr
d&-d&-d&d&0&-d&d&2i+di&-2i-di&0\cr
d&-d&-d&d&0&2d&-2d&0&0&0
\end{pmatrix}\nonumber\\
T&=&\hbox{diag}(1,-1,1,-1,1,e^{\pi i/3}, e^{4\pi i/3},
e^{5\pi i/6},e^{5\pi i/6},-i). \nonumber
\end{eqnarray}
}

\subsection{Verlinde matrices of the quantum $E_6$ model}

We computed numerically the Verlinde matrices of the quantum
double of $E_6$: $N_0$ is the identity matrix, 
and the others can be written as the 
following quadratic expressions
\begin{eqnarray*}
N_1&=&\chi_0\chi_1^\ast+\chi_1\chi_0^\ast+\chi_2\chi_3^\ast+\chi_3\chi_2^\ast+|\chi_4|^2+
\chi_5\chi_6^\ast+\chi_6\chi_5^\ast+\chi_7\chi_8^\ast+
\chi_8\chi_7^\ast+|\chi_9|^2,\\
N_2&=&\chi_0\chi_2^\ast+\chi_2\chi_0^\ast+\chi_1\chi_3^\ast+\chi_3\chi_1^\ast+
|\chi_2+\chi_3|^2+(\chi_2+\chi_3)(\chi_5+\chi_6)^\ast\\ &&+
(\chi_5+\chi_6)(\chi_2+\chi_3)^\ast 
+2|\chi_4|^2+\chi_4(\chi_7+\chi_8+\chi_9)^\ast+
(\chi_7+\chi_8+\chi_9)\chi_4^\ast\\
&&+\chi_5\chi_6^\ast+\chi_6\chi_5^\ast
+\chi_7(\chi_8+\chi_9)^\ast+
(\chi_8+\chi_9)\chi_7^\ast+\chi_8\chi_9^\ast+\chi_9\chi_8^\ast,\\
N_3&=&\chi_0\chi_3^\ast+\chi_3\chi_0^\ast+\chi_1\chi_2^\ast+
\chi_2\chi_1^\ast+|\chi_2+\chi_3|^2+
(\chi_2+\chi_3)(\chi_5+\chi_6)^\ast\\
&&+(\chi_5+\chi_6)(\chi_2+\chi_3)^\ast+2|\chi_4|^2+
\chi_4(\chi_7+\chi_8+\chi_9)^\ast+(\chi_7+\chi_8+\chi_9)\chi_4^\ast\\
&&+|\chi_5|^2+|\chi_6|^2+|\chi_7|^2+|\chi_8|^2+
\chi_9(\chi_7+\chi_8)^\ast+(\chi_7+\chi_8)\chi_9^\ast,\\
N_4&=&(\chi_0+\chi_1+2\chi_2+2\chi_3)\chi_4^\ast+\chi_4(\chi_0+\chi_1+
2\chi_2+2\chi_3)^\ast+\chi_4(\chi_5+\chi_6)^\ast\\
&&+(\chi_5+\chi_6)\chi_4^\ast+(\chi_2+\chi_3+\chi_5+\chi_6)
(\chi_7+\chi_8+\chi_9)^\ast\\
&&+(\chi_7+\chi_8+\chi_9)(\chi_2+\chi_3+\chi_5+\chi_6)^\ast,\\
N_5&=&(\chi_0+\chi_3+\chi_5)\chi_5^\ast+\chi_5(\chi_0+\chi_3)^\ast+
|\chi_2+\chi_3|^2\\
&&+(\chi_4+\chi_7+\chi_8+\chi_9)\chi_4^\ast+\chi_4(\chi_7+\chi_8+\chi_9)^\ast+
(\chi_1+\chi_2+\chi_6)\chi_6^\ast\\
&&+\chi_6(\chi_1+\chi_2+\chi_6)^\ast+\chi_7\chi_8^\ast+
\chi_8\chi_7^\ast+|\chi_9|^2,
\end{eqnarray*}
\begin{eqnarray*}
N_6&=&\chi_0\chi_6^\ast+\chi_6\chi_0^\ast+\chi_1\chi_5^\ast+
\chi_5\chi_1^\ast+|\chi_2+\chi_3|^2+\chi_2\chi_5^\ast+
\chi_5\chi_2^\ast+\chi_3\chi_6^\ast\\
&&+\chi_6\chi_3^\ast+(\chi_4+\chi_7+\chi_8+\chi_9)\chi_4^\ast+
\chi_4(\chi_7+\chi_8+\chi_9)^\ast+\chi_5\chi_6^\ast\\
&&+\chi_6\chi_5^\ast+|\chi_7|^2+|\chi_8|^2+|\chi_9|^2,\\
N_7&=&\chi_0\chi_7^\ast+\chi_8\chi_0^\ast+\chi_1\chi_8^\ast+\chi_7\chi_1^\ast+
\chi_2(\chi_4+\chi_8+\chi_9)^\ast+(\chi_4+\chi_7+\chi_9)\chi_2^\ast\\
&&+\chi_3(\chi_4+\chi_7+\chi_9)^\ast+
(\chi_4+\chi_8+\chi_9)\chi_3^\ast+(\chi_5+\chi_6)\chi_4^\ast+
\chi_4(\chi_5+\chi_6)^\ast\\
&&+\chi_5\chi_8^\ast+
\chi_7\chi_5^\ast+\chi_6\chi_7^\ast+\chi_8\chi_6^\ast,\\
N_8&=&\chi_0\chi_8^\ast+\chi_7\chi_0^\ast+\chi_1\chi_7^\ast+\chi_8\chi_1^\ast
+\chi_2(\chi_4+\chi_7+\chi_9)^\ast\\
&&+(\chi_4+\chi_8+\chi_9)\chi_2^\ast+
(\chi_4+\chi_7+\chi_9)\chi_3^\ast
+\chi_3(\chi_4+\chi_8+\chi_9)^\ast\\
&&+(\chi_5+\chi_6)\chi_4^\ast+
\chi_4(\chi_5+\chi_6)^\ast+\chi_5\chi_7^\ast +\chi_7\chi_6^\ast+
\chi_6\chi_8^\ast+
\chi_8\chi_5^\ast,\\
N_9&=&\chi_0\chi_9^\ast+\chi_9\chi_0^\ast+
(\chi_2+\chi_3)(\chi_4+\chi_7+\chi_8)^\ast
+(\chi_4+\chi_7+\chi_8)(\chi_2+\chi_3)^\ast\\ &&+
(\chi_5+\chi_6)(\chi_4+\chi_9)^\ast+
(\chi_4+\chi_9)(\chi_5+\chi_6)^\ast+\chi_1\chi_9^\ast+\chi_9\chi_1^\ast.
\end{eqnarray*}

All the $N_i$'s are symmetric apart from $N_7$ and $N_8$ that are
transpose of each other. The first two matrices $N_0$ and $N_1$
are permutation matrices as the primary fields 0 and 1 are the
simple currents of our present modular data. The quantum $E_6$
modular invariant $Z_2$ is the charge conjugation invariant, i.e.
${[Z_2]}_{\la,\mu}=\delta_{\la,\bar{\mu}}$. The Frobenius-Schur
indicator FS$_\la=1$ for all $\la\in\cA$ except for $\la=7,8$
where it is zero.

\subsection{Modular invariants of the quantum $E_6$ subfactor(s)}

The dimension of the commutant $\{S,T\}^\prime$ is four and
spanned by the modular invariants which are exactly four and were
computed numerically using \cite[\eerf{1.3}]{BER2}. 
We obtain $Z_1=N_0$,
and the others are given by the 
following quadratic expressions:
\begin{eqnarray*}
Z_2&=&|\ch_0|^2+|\ch_1|^2+|\ch_2|^2+|\ch_3|^2+|\ch_4|^2+|\ch_5|^2+|\ch_6|^2+
|\ch_9|^2+\ch_7\ch_8^\ast+\ch_8\ch_7^\ast,\\
Z_3&=&|\ch_0+\ch_2|^2+|\ch_1+\ch_3|^2+2|\ch_4|^2,\\
Z_4&=&|\ch_0+\ch_2+\ch_4|^2.
\end{eqnarray*}


\subsection{Subfactors for quantum $E_6$ modular invariants}

Here we consider the Longo-Rehren $N\subset M$ inclusion of the
even system of the $E_6$ Dynkin diagram, and let $\cA$ be the
system of endomorphisms that yield the above modular data. By
\cite[Lemma 3.8]{EP} (see also \cite{R3}) we can choose an endomorphism 
$\la_1$ from the sector $[\la_1]$ such that $\la_1^2=1$ 
since
$T_{\la_1,\la_1}=1$. Therefore $[\theta]=[\la_0]\oplus[\la_1]$ is
a dual sector of the braided subfactor $N\subset N\rtimes\bbZ_2$. The 
irreducible decomposition of the $N$-$M$ sectors is as follows: 
$[\iota\la_0]=[\iota\la_1],$
$[\iota\la_2]=[\iota\la_3],$
$[\iota\la_4]=[a_1]\oplus[a_2],$
$[\iota\la_5]=[\iota\la_6],$
$[\iota\la_7]=[\iota\la_8],$ and 
$[\iota\la_9]=[a_3]\oplus[a_4]$. 
Hence gives rise to a trace 8 modular invariant. So $Z_{N\subset
N\rtimes\bbZ_2}=Z_2$ the permutation invariant. Second the
LR$(E_6)$ dual canonical sector is
$[\theta]=[\la_0]\oplus[\la_2]\oplus[\la_4]$. The irreducible
decomposition of the $N$-$M$ sectors are as follows:
$[\iota\la_0], [\iota\la_1], [\iota\la_2]=[\iota\la_0]\oplus
[\iota\la_5]$, $[\iota\la_3]=[\iota\la_1]\oplus[\iota\la_5]$,
$[\iota\la_4]=[\iota\la_0]\oplus[\iota\la_1]\oplus[\iota\la_5]$,
$[\iota\la_5]=[\iota\la_6]=[\iota\la_7]=[\iota\la_8]=[\iota\la_9]$.
Therefore LR$(E_6)$ produces a trace 3 modular invariant. Hence
$Z_{\mathrm{LR}(E_6)}=Z_4.$ By \cite[Sect.\ 5]{BEK3} the notation
of Fig.\ 1 on page 28, $[(0,0)]\oplus[(6,0)]$ is the
$[\la_0]\oplus[\la_2]$ in our current notation. In \cite[Example
5.1]{BEK3} instead of taking the three even vertices of $E_6$ one
can take the two extreme endpoints (which are even and discard
the internal point) following Izumi's Galois correspondence
\cite[Proposition 2.4]{iz1}. In this way, we get an intermediate
subfactor $N\subset P\subset M$ with inclusion map $\iota_P$ such
that $[\theta]=[\la_0]\oplus[\la_2]$ is the sector of the dual
canonical endomorphism of $N\subset P$. Then by the Verlinde
matrices we find that: $[\iota_P\la_0], [\iota_P\la_1]$,
$[\iota_P\la_2]=[\iota_P\la_0]\oplus[\iota_P\la_5]$,
$[\iota_P\la_4]=[\iota_P\la_9]\oplus[a]\oplus[b]$,
$[\iota_P\la_5]=[\iota_P\la_6]$, $[\iota_P\la_7]=[\iota_P\la_8]$,
$[\iota_P\la_9]$, which gives trace 6 modular invariant. Thus the
subfactor $N\subset P$ produces $Z_3$.
Moreover, chiral locality holds true
here. We thus have proven the following.

\begin{theorem}
All the modular invariants of the quantum double of the $E_6$
system are realised by braided subfactors.
\end{theorem}

The fusion rules between the quantum $E_6$ double (sufferable)
modular invariants are as in Table \ref{fusione6}. It is worth
mentioning that as the modular invariants form a basis of the
commutant $\{S,T\}^\prime$ the coefficients exist and are
uniquely determined, but due to a fusion rule algebra phenomenon
those coefficients are non-negative integers \cite[Corollary 3.6]{EP} 
(see also \cite{E1}).

\unitlength 1.4mm \thinlines
\begin{table}[htb]
\begin{center}
\begin{tabular}{|c|c c c c|} \hline
&$Z_1$&$Z_2$&$Z_3$&$Z_4$\\  \hline\hline $Z_1$&
$Z_1$&$Z_2$&$Z_3$&$Z_4$\cr $Z_2$&$Z_2$&$Z_1$ &$Z_3$&$Z_4$\cr
$Z_3$&$Z_3$&$Z_3$&$2Z_3$ &$2Z_4$\cr
$Z_4$&$Z_4$&$Z_4$&$2Z_3$&$3Z_4$\cr \hline
\end{tabular}
\end{center}
\caption{Fusion $Z_aZ_b^t$ of modular invariants} 
\label{fusione6}
\end{table}

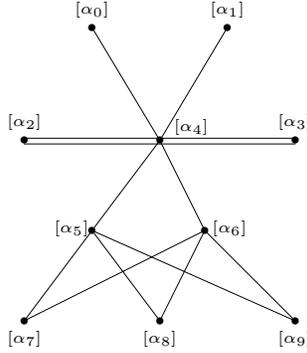
\begin{figure}[htb]
\begin{center}
\unitlength 0.6mm
\begin{picture}(75,65)
\thinlines

\put(20,75){\makebox(0,0){{\tiny $\bullet$}}}
\put(50,75){\makebox(0,0){{\tiny $\bullet$}}}

\put(5,50){\makebox(0,0){{\tiny $\bullet$}}}
\put(35,50){\makebox(0,0){{\tiny $\bullet$}}}
\put(65,50){\makebox(0,0){{\tiny $\bullet$}}}

\put(45,30){\makebox(0,0){{\tiny $\bullet$}}}
\put(20,30){\makebox(0,0){{\tiny $\bullet$}}}

\put(5,10){\makebox(0,0){{\tiny $\bullet$}}}
\put(35,10){\makebox(0,0){{\tiny $\bullet$}}}
\put(65,10){\makebox(0,0){{\tiny $\bullet$}}}

\path(20,75)(35,50)(50,75) \path(5,50.5)(35,50.5)(65,50.5)
\path(5,49.2)(35,49.2)(65,49.2)

\path(35,50)(20,30)(5,10)(45,30)(65,10)(20,30)(35,10)(45,30)
\path(35,50)(45,30)

\put(20,79){\makebox(0,0){{\tiny $[\a_0]$}}}
\put(50,79){\makebox(0,0){{\tiny $[\a_1]$}}}

\put(5,54){\makebox(0,0){{\tiny $[\a_2]$}}}
\put(42,53){\makebox(0,0){{\tiny $[\a_4]$}}}
\put(65,54){\makebox(0,0){{\tiny $[\a_3]$}}}

\put(15.5,30.5){\makebox(0,0){{\tiny $[\a_5]$}}}
\put(50.5,30.5){\makebox(0,0){{\tiny $[\a_6]$}}}

\put(5,6){\makebox(0,0){{\tiny $[\a_7]$}}}
\put(35,6){\makebox(0,0){{\tiny $[\a_8]$}}}
\put(65,6){\makebox(0,0){{\tiny $[\a_9]$}}}

\end{picture}
\end{center}
\caption{$Z_{1}$, fusion graph of $[\a_4^+]$}
\label{e6z1}
\end{figure}

\subsection{Full systems of the quantum $E_6$ modular invariants}

The global index of our system is $\omega=8(1+d+d^2)\sim 89.5692$.

\noindent \subsubsection*{Case $Z_1$}

For the trivial modular invariant $Z_1$ we have
$\cA\simeq\cC^0\simeq\cC.$
 We display the fusion graph of $[\a_4]$ in \fig{e6z1}.

\noindent \subsubsection*{Case $Z_2$}

For the permutation invariant $Z_2$, the full system $\cC$ is
obtained by permuting the sectors of those in the case $Z_1$.

\begin{figure}[htb]
\begin{center}
\unitlength 0.6mm
\begin{picture}(80,110)
\thinlines

\put(0,90){\makebox(0,0){{\tiny $\bullet$}}}
\put(0,75){\makebox(0,0){{\tiny $\bullet$}}}
\put(0,60){\makebox(0,0){{\tiny $\bullet$}}}
\put(0,45){\makebox(0,0){{\tiny $\bullet$}}}

\put(0,90){\arc{2.5}{0}{6.300}} \put(0,60){\arc{3.5}{0}{6.300}}
\put(0,75){\arc{2.5}{0}{6.300}} \put(0,45){\arc{3.5}{0}{6.300}}
\put(0,60){\arc{2.5}{0}{6.300}} \put(0,90){\arc{3.5}{0}{6.300}}
\put(0,45){\arc{2.5}{0}{6.300}} \put(0,45){\arc{3.5}{0}{6.300}}

\put(-25,50){\makebox(0,0){{\tiny $\bullet$}}}
\put(25,50){\makebox(0,0){{\tiny $\bullet$}}}
\put(-18,65){\makebox(0,0){{\tiny $\bullet$}}}
\put(18,65){\makebox(0,0){{\tiny $\bullet$}}}

\put(-25,50){\arc{2.5}{0}{6.300}}
\put(-18,65){\arc{2.5}{0}{6.300}}
\put(25,50){\arc{3.5}{0}{6.300}}
\put(18,65){\arc{3.5}{0}{6.300}}

\put(-10,15){\makebox(0,0){{\tiny $\bullet$}}}
\put(-35,30){\makebox(0,0){{\tiny $\bullet$}}}
\put(7,15){\makebox(0,0){{\tiny $\bullet$}}}
\put(35,30){\makebox(0,0){{\tiny $\bullet$}}}

\path(0,90)(-18,65)(0,75) \path(0,60)(-25,50)(0,45)
\path(-19.4,65)(-26.4,49.9) \path(-18.2,65)(-25.2,50)
\path(-10,15)(18,65)(7,15) \path(-35,30)(25,50)(35,30)

\dottedline(0,90)(18,65) \dottedline(18,65)(0,75)
\dottedline(0,60)(25,50) \dottedline(0,45)(25,50)
\dottedline(24.5,50)(17.5,65) \dottedline(25.4,50)(18.2,65)
\dottedline(-18,65)(-10,15) \dottedline(-18,65)(7,15)
\dottedline(-25,50)(-35,30) \dottedline(-25,50)(35,30)

\put(0,94){\makebox(0,0){{\tiny $[\a_0]$}}}
\put(-25,65){\makebox(0,0){{\tiny $[\a_9^+]$}}}
\put(25,65){\makebox(0,0){{\tiny $[\a_9^-]$}}}
\put(-32,50){\makebox(0,0){{\tiny $[\a_5^+]$}}}
\put(32,50){\makebox(0,0){{\tiny $[\a_5^-]$}}}
\put(0,79){\makebox(0,0){{\tiny $[\a_1]$}}}

\put(0,65){\makebox(0,0){{\tiny $[\a_4^{(1)}]$}}}
\put(0,50){\makebox(0,0){{\tiny $[\a_4^{(2)}]$}}}

\put(-33,27){\makebox(0,0){{\tiny $[(\a_5^+\a_9^-)^{(1)}]$}}}
\put(40,27){\makebox(0,0){{\tiny $[(\a_5^+\a_9^-)^{(2)}]$}}}
\put(-15,5){\makebox(0,0){{\tiny $[(\a_5^+\a_5^-)^{(1)}]$}}}
\put(15,5){\makebox(0,0){{\tiny $[(\a_5^+\a_5^-)^{(2)}]$}}}

\put(110,90){\makebox(0,0){{\tiny $\bullet$}}}
\put(110,75){\makebox(0,0){{\tiny $\bullet$}}}
\put(110,60){\makebox(0,0){{\tiny $\bullet$}}}
\put(110,45){\makebox(0,0){{\tiny $\bullet$}}}

\put(110,90){\arc{2.5}{0}{6.300}}
\put(110,90){\arc{3.5}{0}{6.300}}
\put(110,75){\arc{2.5}{0}{6.300}}
\put(110,75){\arc{3.5}{0}{6.300}}
\put(110,60){\arc{2.5}{0}{6.300}}
\put(110,60){\arc{3.5}{0}{6.300}}
\put(110,45){\arc{2.5}{0}{6.300}}
\put(110,45){\arc{3.5}{0}{6.300}}

\put(85,50){\makebox(0,0){{\tiny $\bullet$}}}
\put(135,50){\makebox(0,0){{\tiny $\bullet$}}}
\put(92,65){\makebox(0,0){{\tiny $\bullet$}}}
\put(128,65){\makebox(0,0){{\tiny $\bullet$}}}

\put(85,50){\arc{2.5}{0}{6.300}} \put(92,65){\arc{2.5}{0}{6.300}}
\put(135,50){\arc{3.5}{0}{6.300}}
\put(128,65){\arc{3.5}{0}{6.300}}

\put(100,15){\makebox(0,0){{\tiny $\bullet$}}}%
\put(90,30){\makebox(0,0){{\tiny $\bullet$}}}
\put(120,15){\makebox(0,0){{\tiny $\bullet$}}}%
\put(133,30){\makebox(0,0){{\tiny $\bullet$}}}

\path(110,75)(92,65)(110,90) \path(110,60)(85,50)(110,45)
\path(90,30)(128,65)(133,30) \path(100,15)(135,50)(120,15)

\dottedline(110,90)(128,65) \dottedline(128,65)(110,75)
\dottedline(135,50)(110,60) \dottedline(135,50)(110,75)
\dottedline(92,65)(90,30) \dottedline(92,65)(133,30)
\dottedline(85,50)(100,15) \dottedline(85,50)(120,15)

\put(89.5,65){\ellipse{5}{3}} \put(88.5,65){\ellipse{8.5}{3}}
\put(82,50){\ellipse{5}{3}} \put(81,50){\ellipse{8.5}{3}}

\path(128,65)(128.9,65.6) \path(129.5,65.8)(130.5,66.1)
\path(129.3,66)(130.5,66.6) \path(128,65)(129,64.5)
\path(129.3,64.2)(130.8,63.5) \path(129.4,64.2)(130.8,64)
\put(133,65){\arc{3}{4.1}{8.5}} \put(133,65){\arc{5}{4}{8.5}}

\path(135,50)(135.9,50.6) \path(136,51)(136.6,51.5)
\path(136,50.9)(136.6,51)

\path(135,50)(135.9,50.6) \path(136.3,49.2)(137.3,48.5)
\path(136.2,49.2)(137.3,49) \put(139,50){\arc{3}{4.1}{8.5}}
\put(139,50){\arc{5}{4}{8.5}}

\put(110,94){\makebox(0,0){{\tiny $[\a_0]$}}}
\put(90,70){\makebox(0,0){{\tiny $[\a_9^+]$}}}
\put(135,70){\makebox(0,0){{\tiny $[\a_9^-]$}}}
\put(77,55){\makebox(0,0){{\tiny $[\a_5^+]$}}}
\put(142,55){\makebox(0,0){{\tiny $[\a_5^-]$}}}
\put(110,80){\makebox(0,0){{\tiny $[\a_1]$}}}

\put(110,65){\makebox(0,0){{\tiny $[\a_4^{(1)}]$}}}
\put(110,40){\makebox(0,0){{\tiny $[\a_4^{(2)}]$}}}

\put(80,27){\makebox(0,0){{\tiny $[(\a_5^+\a_9^-)^{(1)}]$}}}
\put(135,25){\makebox(0,0){{\tiny $[(\a_5^+\a_9^-)^{(2)}]$}}}
\put(93,5){\makebox(0,0){{\tiny $[(\a_5^+\a_5^-)^{(1)}]$}}}
\put(127,5){\makebox(0,0){{\tiny $[(\a_5^+\a_5^-)^{(2)}]$}}}

\end{picture}
\end{center}
\caption{$Z_3$, fusion graphs of $[\a_9^\pm]$ and $[\a_5^\pm]$ where $Z_3^2=2Z_3$}
\label{qse6}
\end{figure}
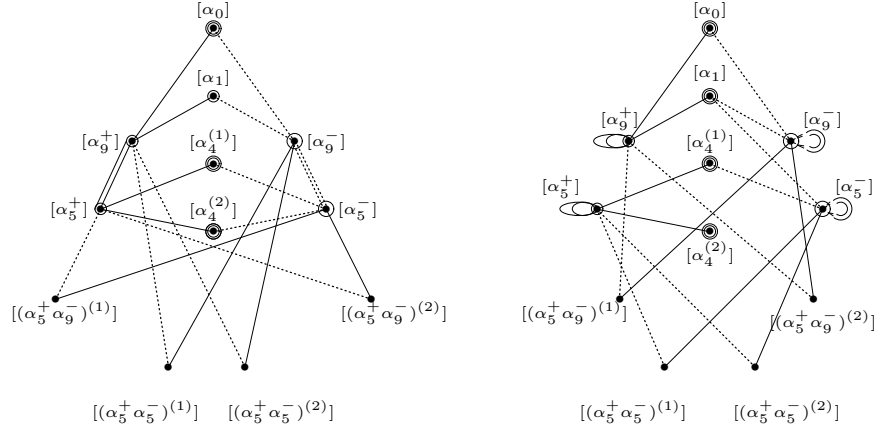

\noindent \subsubsection*{Case $Z_3$}

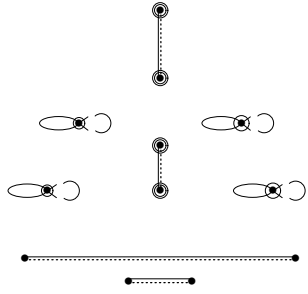
\begin{figure}[htb]
\begin{center}
\unitlength 0.6mm
\begin{picture}(75,65)
\thinlines

\put(30,55){\makebox(0,0){{\tiny $\bullet$}}}
\put(30,40){\makebox(0,0){{\tiny $\bullet$}}}
\put(30,25){\makebox(0,0){{\tiny $\bullet$}}}
\put(30,15){\makebox(0,0){{\tiny $\bullet$}}}

\put(30,55){\arc{2.5}{0}{6.300}} \put(30,55){\arc{3.5}{0}{6.300}}
\put(30,40){\arc{2.5}{0}{6.300}} \put(30,40){\arc{3.5}{0}{6.300}}
\put(30,25){\arc{2.5}{0}{6.300}} \put(30,25){\arc{3.5}{0}{6.300}}
\put(30,15){\arc{2.5}{0}{6.300}} \put(30,15){\arc{3.5}{0}{6.300}}

\put(5,15){\makebox(0,0){{\tiny $\bullet$}}}
\put(55,15){\makebox(0,0){{\tiny $\bullet$}}}
\put(12,30){\makebox(0,0){{\tiny $\bullet$}}}
\put(48,30){\makebox(0,0){{\tiny $\bullet$}}}

\put(5,15){\arc{2.5}{0}{6.300}} \put(12,30){\arc{2.5}{0}{6.300}}
\put(55,15){\arc{3.5}{0}{6.300}} \put(48,30){\arc{3.5}{0}{6.300}}

\put(23,-5){\makebox(0,0){{\tiny $\bullet$}}}
\put(0,0){\makebox(0,0){{\tiny $\bullet$}}}
\put(37,-5){\makebox(0,0){{\tiny $\bullet$}}}
\put(60,0){\makebox(0,0){{\tiny $\bullet$}}}

\path(29.6,55)(29.6,40) \dottedline(30.2,55)(30.2,40)
\path(29.6,25)(29.6,15) \dottedline(30.2,25)(30.2,15)
\path(0,0.4)(60,0.4) \dottedline(0,-0.4)(60,-0.4)
\path(23,-4.6)(37,-4.6) \dottedline(23,-5.4)(37,-5.4)

\put(7.5,30){\ellipse{8.5}{3}} \put(43.5,30){\ellipse{8.5}{3}}
\put(0.5,15){\ellipse{8.5}{3}} \put(50.5,15){\ellipse{8.5}{3}}

\put(17,30){\arc{4.3}{4}{8.5}} \put(10,15){\arc{4.3}{4}{8.5}}
\put(53,30){\arc{4.3}{4}{8.5}} \put(60,15){\arc{4.3}{4}{8.5}}

\path(12,30)(14,31.4) \path(12,30)(14,28.5) \path(5,15)(7,16.4)
\path(5,15)(7,13.5)

\path(48,30)(50,31.4) \path(48,30)(50,28.5) \path(55,15)(57,16.4)
\path(55,15)(57,13.5)

\end{picture}
\end{center}
\caption{$Z_3$, fusion graph of $[\a_1^\pm]$ where $Z_3^2=2Z_3$} \label{qse6a1}
\end{figure}

For the sufferable quantum $E_6$ modular invariant $Z_3$ we have by \cite{BEK2}:
$$\omega_\pm=8(1+d+d^2)/(2+d)\sim 18.9282,\quad \omega_0=4.$$
Computing using the Verlinde matrices we find:
$\lan\a_0,\a_0\ran=\lan\a_1^\pm,\a_1^\pm\ran$
$=\lan\a_5^\pm,\a_5^\pm\ran$ $=\lan\a_6^\pm,\a_6^\pm\ran$
$=\lan\a_7^\pm,\a_7^\pm\ran=\lan\a_8^\pm,\a_8^\pm\ran=1,$
$\lan\a_2^\pm,\a_2^\pm\ran=\lan\a_3^\pm,\a_3^\pm\ran=2,$
$\lan\a_4^\pm,\a_4^\pm\ran=3,$
$\lan\a_0,\a_i^\pm\ran=\delta_{0,i}+\delta_{2,i},\quad$
$\lan\a_1^\pm,\a_i^\pm\ran= \delta_{1,i}+\delta_{3,i}$
$\lan\a_2^\pm,\a_3^\pm\ran=\lan\a_2^\pm,\a_5^\pm\ran$
$=\lan\a_2^\pm,\a_6^\pm\ran=\lan\a_3^\pm,\a_5^\pm\ran$
$=\lan\a_3^\pm,\a_6^\pm\ran=1,$ $\lan\a_4^\pm,\a_7^\pm\ran=$
$\lan\a_4^\pm,\a_8^\pm\ran=$ $\lan\a_4^\pm,\a_9^\pm\ran=1,$
$\lan\a_5^\pm,\a_6^\pm\ran=$ $\lan\a_7^\pm,\a_8^\pm\ran=$
$\lan\a_7^\pm,\a_9^\pm\ran=$ $\lan\a_8^\pm,\a_9^\pm\ran=1$,
 and the others vanish.
We then conclude that $[\a_0]$ and $[\a_1^\pm]$ are irreducible, and
$[\a_2^\pm]=[\a_0]\oplus[\a_5^\pm],\quad$
$[\a_3^\pm]=[\a_1^\pm]\oplus[\a_5^\pm],$ $[\a_4^\pm]=$
$[\a_9^\pm]\oplus[\a_4^{\pm(1)}]\oplus[\a_4^{\pm(2)}],$
$[\a_5^\pm]=[\a_6^\pm],\quad [\a_7^\pm]$ $=[\a_8^\pm]=[\a_9^\pm]$.
Using also the matrix $Z$=$Z_3$ and noting that
$[\a_4^{\pm(i)}]=[\a_4^{\mp(i)}]$ for $i=1,2$ we conclude that
$\cC^0=\{\a_0,\a_1,$ $\a_4^{(1)},\a_4^{(2)}\}$, and 
$\cC^\pm=\cC^0\cup\{\a_5^\pm,\a_9^\pm\}$. Moreover since
$$\lan\a_4^+\a_1^-,\a_4^+\ran=
N_{4,4}^\xi Z_{1,\xi}=3=\lan\a_4^+\a_1^-,\a_4^+\a_1^-\ran$$ we
have the sectors of $\cC^0$ equal $\bbZ_2\times\bbZ_2$ with
$[\a_4^{(i)}] [\a_4^{(i)}]=[\mathrm{id}]$. We need to find four
more irreducible sectors in the full system $\cC$ because
$\#\cC=\Tr(ZZ^t)=12$. We compute and find that
\begin{eqnarray}
\lan \a_5^+\a_5^-,\a_5^+\a_5^-\ran=\sum N_{5,5}^\xi N_{5,5}^\eta
Z_{\xi,\eta}=2 \nonumber
\end{eqnarray}
and the dimensions of the intertwiner spaces between
$\a_5^+\a_5^-$ and any sector in both chiral systems $\cC^\pm$
vanish. Hence $[\a_5^+\a_5^-]$ decomposes into two new
irreducible sectors $[(\a_5^+\a_5^-)^{(1)}]$ and
$[(\a_5^+\a_5^-)^{(2)}]$. The other two sectors are similarly
obtained from the decomposition of $\a_9^+\a_5^-$ further noting
that $\lan\a_5^+\a_5^-,\a_9^+\a_5^-\ran=0$. We denote these new
irreducible sectors by $[(\a_9^+\a_5^-)^{(1)}]$ and
$[(\a_9^+\a_5^-)^{(2)}]$. We have the fusion graphs of both
$[\a_9^+]$ and $[\a_9^-]$ displayed on the LHS of \fig{qse6}. We use
straight lines for the fusion graph of $[\a_9^+]$ and dashed
lines for that of $[\a_9^-]$. On the RHS of \fig{qse6} the fusion
graph of $[\a_5^+]$ and $[\a_5^-]$ is displayed. 
We display the
fusion graph of $[\a_1^+]$ and $[\a_1^-]$ in \fig{qse6a1}.
When $\theta$ is the local dual canonical endomorphism 
$\la_0\oplus \la_2$,
let $\gamma$
denote the canonical sector of $N\subset M$. Since
$\lan\a_5^+\a_5^-,\gamma\ran=\lan\a_5^+,\a_5^+\ran=1$ we conclude
from \cite[Corollary 3.19]{BE3} that
$[\gamma]=[\a_0]\oplus[(\a_5^+\a_5^-)^{(1)}]$. 
The full system $\cC$ decomposes as two sheets.
The first is the chiral system $\cC^+$ and the second sheet
comprises the irreducible components of $\cC^+\alpha_5^-$
in accordance with the decomposition of $Z_3^2 = 2Z_3$ 
and $[\gamma]=[\a_0]\oplus[(\a_5^+\a_5^-)^{(1)}]$
and \cite[Corollary 3.6]{EP} and \cite{E1}.

\noindent \subsubsection*{Case $Z_4$}

The global indices of the chiral systems $\cC^\pm$ and $\cC^0$ are by \cite{BEK2}:
$$\omega_\pm=8(1+d+d^2)/(4+2d)=2+d^2\sim 9.4641,\quad 
\omega_0=\omega_\pm^2/\omega=1$$
respectively. There are $9=\Tr(Z_4Z_4^t)$ $M$-$M$ irreducible
sectors in the full commutative chiral system $\cC$.
Note that $\cC^0=\{\a_0\}$ since
$\omega_0=1$. We next compute the chiral systems $\cC^\pm$.
We get (using the streamlined notation
$\a_{\la_i}^\pm=:\a_i^\pm$): $\lan\a_0^\pm,\a_0^\pm\ran$
$=\lan\a_1^\pm,\a_1^\pm\ran=$ $\lan\a_4^\pm,\a_4^\pm\ran=$
$\lan\a_5^\pm,\a_5^\pm\ran=$ $\lan\a_6^\pm,\a_6^\pm\ran$
$=\lan\a_7^\pm,\a_7^\pm\ran=$ $\lan\a_8^\pm,\a_8^\pm\ran$
$=\lan\a_9^\pm,\a_9^\pm\ran=1,$ $\lan\a_2^\pm,\a_2^\pm\ran=$
$\lan\a_3^\pm,\a_3^\pm\ran=2,$ $\lan\a_4^\pm,\a_4^\pm\ran=3,$
$\lan\a_0^\pm,\a_i^\pm\ran=\delta_{0,i}+\delta_{2,i}+\delta_{4,i},\quad$
$\lan\a_1^\pm,\a_i^\pm\ran=\delta_{1,i}+\delta_{3,i}+\delta_{4,i},\quad$
$\lan\a_2^\pm,\a_3^\pm\ran=\lan\a_2^\pm,\a_5^\pm\ran=1,$
$\lan\a_2^\pm,\a_4^\pm\ran=2,$ $\lan\a_3^\pm,\a_4^\pm\ran=2,$
$\lan\a_3^\pm,\a_5^\pm\ran=\lan\a_3^\pm,\a_6^\pm\ran=1,$
$\lan\a_4^\pm,\a_5^\pm\ran=\lan\a_4^\pm,\a_6^\pm\ran=$
$\lan\a_5^\pm,\a_6^\pm\ran=\lan\a_5^\pm,\a_7^\pm\ran$
$=\lan\a_5^\pm,\a_8^\pm\ran$ $=\lan\a_5^\pm,\a_9^\pm\ran=1.$

Therefore $\cC^\pm=\{\a_0,\a_1^\pm,\a_5^\pm\}$ and as sectors
$[\a_2^\pm]=$ $[\a_0]\oplus[\a_5^\pm],$ $[\a_3^\pm]=$
$[\a_1^\pm]\oplus[\a_5^\pm],$ $[\a_4^\pm]=$
$[\a_0]\oplus[\a_1^\pm]\oplus[\a_5^\pm],$ $[\a_5^\pm]=$
$[\a_6^\pm]=$ $[\a_7^\pm]=$ $[\a_8^\pm]=$ $[\a_9^\pm]$.

\begin{figure}[htb]
\begin{center}
\unitlength 0.6mm
\begin{picture}(75,110)
\thinlines

\put(0,0){\makebox(0,0){{\tiny $\bullet$}}}
\put(0,15){\makebox(0,0){{\tiny $\bullet$}}}
\put(0,90){\makebox(0,0){{\tiny $\bullet$}}}

\put(-20,15){\makebox(0,0){{\tiny $\bullet$}}}
\put(-20,15){\makebox(0,0){{\tiny $\bullet$}}}
\put(25,35){\makebox(0,0){{\tiny $\bullet$}}}
\put(-25,35){\makebox(0,0){{\tiny $\bullet$}}}

\put(-25,65){\makebox(0,0){{\tiny $\bullet$}}}
\put(25,65){\makebox(0,0){{\tiny $\bullet$}}}

\path(-25,65)(-25,35)(0,90) \path(0,15)(-20,15)(25,65)
\path(20,15)(0,0)(25,35) \dottedline(0,90)(25,35)
\dottedline(25,35)(25,65) \dottedline(-25,65)(20,15)
\dottedline(20,15)(0,15) \dottedline(-25,35)(0,0)
\dottedline(0,0)(-20,15)

\put(-25,65){\arc{2.5}{0}{6.300}}
\put(-25,35){\arc{2.5}{0}{6.300}} \put(0,90){\arc{2.5}{0}{6.300}}
\put(0,90){\arc{3.5}{0}{6.300}} \put(25,65){\arc{3.5}{0}{6.300}}
\put(25,35){\arc{3.5}{0}{6.300}}

\put(-27.5,35){\ellipse{5}{3}} \put(-18.5,15){\ellipse{5}{3}}
\put(-3.5,0){\ellipse{5}{3}}

\put(-28.7,35){\ellipse{8.5}{3}} \put(-23.7,15){\ellipse{8.5}{3}}
\put(-4,0){\ellipse{8.5}{3}}

\path(2,0.5)(0,0)(2,-1.7) \path(2.5,0.5)(3.5,0.5)
\path(2.5,-1.9)(3.5,-2) \put(5,-0.5){\arc{3}{4.1}{8.5}}
\put(5,-0.5){\arc{5}{4}{8.5}}

\path(2,0.9)(0,0)(2,-2.3) \path(2.4,1.1)(3,1.3)
\path(2.5,-2.4)(2.9,-2.5)
\path(22,16)(20,15)(22,14) \path(22.5,16)(23.5,16)
\path(22.5,14)(23.5,14) \put(25,15){\arc{3}{4.1}{8.5}}
\put(25,15){\arc{5}{4}{8.5}}

\path(22,16.4)(20,15)(22,13.5) \path(22.4,16.6)(22.9,16.7)
\path(22.5,13.3)(22.9,13.2)

\path(27,36)(25,35)(27,34) \path(27.5,36)(28.5,36)
\path(27.5,34)(28.5,34) \put(30,35){\arc{3}{4.1}{8.5}}
\put(30,35){\arc{5}{4}{8.5}}

\path(27,36.4)(25,35)(27,33.5) \path(27.4,36.6)(27.9,36.8)
\path(27.5,33.2)(27.9,33.1)

\put(0,-5.5){\makebox(0,0){{\tiny $[\a_5^+\a_5^-]$}}}
\put(0,18.5){\makebox(0,0){{\tiny $[\a_1^+\a_1^-]$}}}
\put(0,95){\makebox(0,0){{\tiny $[\hbox{id}]$}}}

\put(26,10){\makebox(0,0){{\tiny $[\a_1^+\a_5^-]$}}}
\put(-25,10.5){\makebox(0,0){{\tiny $[\a_5^+\a_1^-]$}}}
\put(31,41){\makebox(0,0){{\tiny $[\a_5^-]$}}}
\put(-32,41){\makebox(0,0){{\tiny $[\a_5^+]$}}}

\put(-25,70){\makebox(0,0){{\tiny $[\a_1^+]$}}}
\put(25,70){\makebox(0,0){{\tiny $[\a_1^-]$}}}


\put(90,0){\makebox(0,0){{\tiny $\bullet$}}}
\put(90,15){\makebox(0,0){{\tiny $\bullet$}}}
\put(90,90){\makebox(0,0){{\tiny $\bullet$}}}
\put(110,15){\makebox(0,0){{\tiny $\bullet$}}}
\put(70,15){\makebox(0,0){{\tiny $\bullet$}}}
\put(115,35){\makebox(0,0){{\tiny $\bullet$}}}
\put(65,35){\makebox(0,0){{\tiny $\bullet$}}}
\put(65,65){\makebox(0,0){{\tiny $\bullet$}}}
\put(115,65){\makebox(0,0){{\tiny $\bullet$}}}

\put(65,65){\arc{2.5}{0}{6.300}} \put(65,35){\arc{2.5}{0}{6.300}}
\put(90,90){\arc{2.5}{0}{6.300}} \put(90,90){\arc{3.5}{0}{6.300}}
\put(115,65){\arc{3.5}{0}{6.300}}
\put(115,35){\arc{3.5}{0}{6.300}}

\path(90,90)(65,65) \path(115,65)(90,15) \path(115,35)(110,15)
\dottedline(90,90)(115,65) \dottedline(65,65)(90,15)
\dottedline(65,35)(70,15)

\put(60.5,35){\ellipse{8.5}{3}} \put(65.5,15){\ellipse{8.5}{3}}
\put(85.5,0){\ellipse{8.5}{3}}

\path(92,0.5)(90,0)(92,-1.7) \path(92.5,0.5)(93.5,0.5)
\path(92.5,-1.9)(93.5,-2) \put(95,-0.5){\arc{3}{4.1}{8.5}}

\path(112,16)(110,15)(112,14) \path(112.5,16)(113.5,16)
\path(112.5,14)(113.5,14) \put(115,15){\arc{3}{4.1}{8.5}}

\path(117,36)(115,35)(117,34) \path(117.5,36)(118.5,36)
\path(117.5,34)(118.5,34) \put(120,35){\arc{3}{4.1}{8.5}}

\put(90,-5.5){\makebox(0,0){{\tiny $[\a_5^+\a_5^-]$}}}
\put(90,10.5){\makebox(0,0){{\tiny $[\a_1^+\a_1^-]$}}}
\put(90,95){\makebox(0,0){{\tiny $[\hbox{id}]$}}}

\put(116,10){\makebox(0,0){{\tiny $[\a_1^+\a_5^-]$}}}
\put(65,10.5){\makebox(0,0){{\tiny $[\a_5^+\a_1^-]$}}}
\put(121,41){\makebox(0,0){{\tiny $[\a_5^-]$}}}
\put(58,41){\makebox(0,0){{\tiny $[\a_5^+]$}}}

\put(65,70){\makebox(0,0){{\tiny $[\a_1^+]$}}}
\put(115,70){\makebox(0,0){{\tiny $[\a_1^-]$}}}

\end{picture}
\end{center}
\caption{$Z_{4}$, fusion graphs of $[\a_5^\pm]$ and $[\a_1^\pm]$ where $Z_4^2=3Z_4$}
\label{e6z4}
\end{figure}
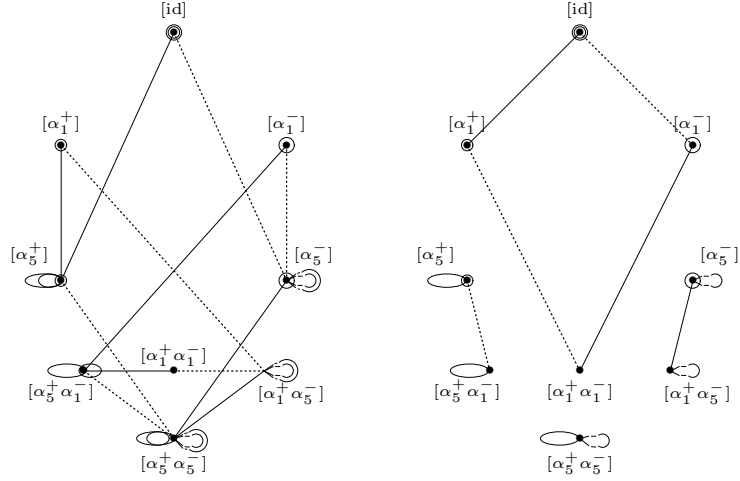

The fusion rules are as the well known fusion rules of the even
vertices of the graph $E_6$ \cite{BEK3}: $[\a_5^\pm\a_5^\pm]=$
$[\a_0]\oplus[\a_1^\pm]\oplus 2[\a_5^\pm]$,
$[\a_1^\pm][\a_5^\pm]=[\a_5^\pm]$, $[\a_1^\pm][\a_1^\pm]=[\a_0].$
We find that $\cC=\cC^+\times\cC^-$ and remark that
$d_{\a_5^\pm}=1+\sqrt{3}$. We display the fusion the all $M$-$M$
system together with the fusion graphs of both $[\a_5^-]$ and
$[\a_5^+]$ in the LHS \fig{e6z4}. In this figure, we use straight
lines for the fusion graph of $[\a_5^+]$ whereas dashed lines for
that of $[\a_5^-]$. We also encircled the $\cC^+$-chiral
sectors with small circles and with larger circles those for the
$\cC^-$-chiral system. The RHS of \fig{e6z4} we display
the fusion graphs of $[\a_1^\pm]$. 
If
$\theta$ is the local dual canonical endomorphism $\la_0\oplus
\la_2 \oplus \la_4$,
let $\gamma$ be the corresponding canonical
endomorphism of $N\subset M$. As an application of
\cite[Corollary 3.19]{BE3}, we can compute that
$[\gamma]=[\a_0]\oplus[\a_1^+\a_1^-]\oplus[\a_5^+\a_5^-]$.
In accordance with $Z_4^2 = 3Z_4$ and 
\cite[Corollary 3.6]{EP} the full system
$\cC$
decomposes as three sheets $\cC^+$, $\cC^+\alpha_1^-$ and
$\cC^+\alpha_5^-$.

\newcommand\bitem[2]{\bibitem{#1}{#2}}

\def\aam              {Acta Appl.\ Math. }
\def\aip              {Ann.\ Inst.\ H.\ Poincar\'e (Phys.\ Th\'eor.) }
\def\cmp              {Comm.\ Math.\ Phys. }
\def\duke             {Duke Math.\ J. }
\def\ijm              {Intern.\ J. Math. }
\def\jfa              {J.\ Funct.\ Anal. }
\def\jmp              {J.\ Math.\ Phys. }
\def\lmp              {Lett.\ Math.\ Phys. }
\def\rmp              {Rev.\ Math.\ Phys. }
\def\inv              {Invent.\ Math. }
\def\mpl              {Mod.\ Phys.\ Lett. }
\def\nup              {Nuclear Phys. }
\def\nupb             {Nuclear Phys.\ {\bf B}}
\def\nupp             {Nuclear Phys.\ (Proc.\ Suppl.) }
\def\adma             {Adv.\ Math. }
\def\physa            {Physica \textbf{A} }
\def\ijmp             {Int.\ J.\ Mod.\ Phys. }
\def\jp               {J.\ Phys. }
\def\fdp              {Fortschr.\ Phys. }
\def\pl               {Phys.\ Lett.}
\def\rims             {Publ.\ Res.\ Inst.\ Sci. }
\def\RMP              {Rev.\ Math.\ Phys.}


\bibliographystyle{amsalpha}

\end{document}